\documentclass[12pt,twoside,a4paper]{article}

\usepackage{amssymb}
\usepackage{rotating}
\usepackage{amssymb}
\usepackage{amsmath,amscd}
\usepackage{fullpage}
\usepackage[english]{babel}
\usepackage{pst-node,pstricks,multido,pst-plot,pst-text,pst-3d}
\usepackage{graphicx}
\usepackage{amsfonts}
\usepackage[latin1]{inputenc}
\newtheorem{example}{Example}[section]
\newtheorem{note}[example]{Note}

\newtheorem{corollary}[example]{Corollary}

\newtheorem{proposition}[example]{Proposition}

\newtheorem{lemma}[example]{Lemma}

\def\S{{\mathfrak  S}}
\def\cal#1{{\mathfrak #1}}

\def\<{\langle}
\def\>{\rangle}

\def\C{{\mathbb C}}

\def\Y{{\mathbb Y}}

\def\ashuff#1#2#3{
\kern 1pt \vrule height#1 \overline{\vrule height#3 width 0pt
\hskip#2} \rule{.3pt}{#1}\overline{\vrule height#3 width 0pt
\hskip#2} \rule{.3pt}{#1} \kern 1pt }

\def\X{{\mathbb X}}\def\Y{{\mathbb Y}}

\title{How to compute Selberg-like integrals?}
\author{M. Deneufch\^atel\footnote{LIPN, UMR 7030, Universit\'e Paris 13 - CNRS; 99, Avenue J.-B. Cl\'ement, 93430 Villetaneuse, France}}

\begin{document}
\maketitle
\begin{abstract}

In this paper, we describe a general method for computing Selberg-like integrals based on a formula, due to Kaneko, for Selberg-Jack integrals. The general principle consists in expanding the integrand \emph{w.r.t.} the Jack basis, which is obtained by a Gram-Schmidt orthogonalization process.

The resulting algorithm is not very efficient because of this decomposition. But for special cases, the coefficients admit a closed form. As an example, we study the case of the power-sums since for which the coefficients are obtained by manipulating generating series by means of transformations of alphabets. Furthermore, we prove that the integral is a rational function in the number of variables which allows us to study asymptotics. 
As an application, we investigate the asymptotic behavior when the integrand involves Jack polynomials and power sums.
\end{abstract}

\section{Introduction}
This paper is devoted to Selberg-like integrals, that is:
\begin{equation}
\label{intgen}
 \langle f \rangle_{a,b,\kappa}^{N} = \frac{1}{N!}  \int_{\left[ 0 , 1 \right]^{N}} f(x_{1} , \dots , x_{N}) \prod_{i<j} (x_{i}-x_{j})^{2 \kappa} \prod_{i=1}^{N} x_{i}^{a-1} \left( 1 - x_{i} \right)^{b-1} dx_{i},
\end{equation}
where $f$ is a multivariate polynomial.\\
For a review on Selberg integrals, see \cite{oai:arXiv.org:0710.3981}.

This integral appears in the physical problem of \emph{quantum transport}. A link has been found between the properties of some large random matrices and the fluctuations of the electronic conductance in certain disordered systems (see \cite{RevModPhys.69.731}, \cite{Guhr:1997ve}, \cite{PhysRevLett.59.2475} or \cite{khoruzhenko-2009-80}). In the case of a system constituted by a phase coherent disordered region connected to two reservoirs of electrons by two ideal leads, it is possible to approximate the \emph{scattering matrix} (the matrix that links the wave functions of the incoming and outgoing electrons) by a random matrix $S$ which belongs to one of Dyson's circular ensembles. As usual, the eigenvalues of this random matrix are of great interest, and some experimentaly measurable variables can be computed with them. The contraints of the physical problem restrain the admissible form for the joint probability density of the eigenvalues of $S$ and this leads to the computation of integrals of the form $\langle f \rangle_{a,b,\kappa}^{N}$. The properties of large physical systems are related to the limit $N \rightarrow \infty$.

In \cite{oai:arXiv.org:0912.1228}, Luque and Vivo described a method for computing the asymptotics when $N\rightarrow \infty$ of $\displaystyle{\langle x^{k} \rangle_{a,b,\kappa}^{N}}$. Their method involves different basis of symmetric functions. Surprisingly, the result has an interesting combinatorial interpretation in terms of Dyck paths.\\
The main tool of the algorithm presented in \cite{oai:arXiv.org:0912.1228} is the expansion of power sums 
\begin{equation}
\nonumber
p_{k}(\mathbb{X}) := x_{1}^{k} + \dots + x_{N}^{k} 
\end{equation}
in the Jack basis $P_{\lambda}^{\left( \frac{1}{\kappa} \right)}$. The coefficients arising in the decomposition
\begin{equation}
   p_{k} = \sum_{\lambda \vdash k} \alpha_{\lambda,k} P_{\lambda}^{\left( \frac{1}{\kappa} \right)}
\end{equation}
do not depend on the number of variables $N$. Now, by the well-known formula due to Kaneko (see \cite{Kadell},\cite{Kaneko},\cite{Macdonald95}):
\begin{equation}
\label{intMacdo}
   \langle P_{\lambda}^{\left( \frac{1}{\kappa} \right)} \rangle_{a,b,\kappa}^{N} = \prod_{i<j} \frac{\Gamma(\lambda_{i}-\lambda_{j}+\kappa(j-i+1))}{\Gamma(\lambda_{i}-\lambda_{j}+\kappa(j-i))} \prod_{i=1}^{N} \frac{\Gamma(\lambda_{i} + a + \kappa(N-i)) \Gamma(b+\kappa(N-i))}{\Gamma(\lambda_{i}+a+b+\kappa(2N-i-1))},
\end{equation}
Luque and Vivo obtained the exact value of $\displaystyle{\frac{\langle x^{k} \rangle_{a,b,\kappa}^{N}}{\langle 1 \rangle_{a,b,\kappa}^{N}}}$ with its asymptotic behavior. But they did not explain how to obtain the coefficients $\alpha_{\lambda,k}$. 

In general, Jack polynomials are obtained by orthogonalizing the Schur basis \emph{w.r.t.} a deformation of the usual scalar product. But this does not provide an efficient algorithm. We show how to improve it by computing directly the value of the coefficient $\alpha_{\lambda,k}$.

This paper is organised as follows. In the next section, we recall some basic properties of Jack and Macdonald polynomials. Section \ref{exactcomputation} contains the computation of (\ref{intgen}) for Jack polynomials. In the last section, we give the exact value of the Selberg - power sum integral and present two asymptotic relations suggested by numerical results.

\section{Jack and Macdonald polynomials}
In this section, we recall some definitions and properties of symmetric functions and Jack and Macdonald polynomials that will be used essentially in section \ref{spec}.
\subsection{Symmetric functions, Cauchy Kernel and $\lambda$-ring operations}
\label{cauchykern}
It is well known (see \emph{e.g.} \cite{Lascoux03}) that the algebra of symmetric functions is endowed with a $\lambda$-ring structure. As usual, we will denote by $\sigma_z(\mathbb{X})$ the Cauchy function, that is the generating function of the complete functions $S_i$ over the alphabet $\X=\{x_1,\dots,x_N\}$:
\begin{equation}
\sigma_{z}(\mathbb{X}) = \sum_{i} S_{i}(\mathbb{X}) z^{i} = \prod_{x \in \mathbb{X}} \frac{1}{1-xz}.
\end{equation}
Remark that $\sigma_{z}$ and the exponential generating series of the power sums are related as follows:
\begin{equation}
   \sigma_{z}(\mathbb{X}) = \exp \left( \sum_{n \geq 1} \frac{p_{n}}{n} z^{n} \right).
\end{equation}
The sum of two alphabets $\mathbb{X}$ and $\mathbb{Y}$ is defined by
\begin{equation}
\label{somme}
   \sigma_{z}(\mathbb{X} + \mathbb{Y}) = \sigma_{z}(\mathbb{X}) \sigma_{z}(\mathbb{Y}) = \sum_{i} S_{i}(\mathbb{X}+\mathbb{Y}) z^{i}.
\end{equation}
The multiplication of an alphabet by a constant is defined by extending to any complex number $u$ the following property: $\displaystyle{\sigma_{z}(2 \mathbb{X}) = \sigma_{z}^{2}(\mathbb{X})}$ which is obtained with $\mathbb{X} = \mathbb{Y}$ in (\ref{somme}). Thus, 
\begin{equation}
\sigma_{z}(u \mathbb{X}) = \sigma_{z}^{u}(\mathbb{X}).
\end{equation}
In particular, if $u=-1$, $\sigma_z(-\X)=\sigma_z(\X)^{-1}$. Hence, multiplicities could appear in the alphabet: the notion of alphabet is not restricted to sets of variables but is, by this way, extended to series. For example, we will use the alphabet $\frac{1-u}{1-t}$, which will be considered as the difference of the alphabets $1+t+\dots+t^n+\dots$ and $u+tu+\dots+t^nu+\dots$. In terms of operations on symmetric functions, this is equivalent to the operation that maps $p_{n}$ on $\displaystyle{\frac{1-u^{n}}{1-t^{n}}}$.\\
Finally, the product of two alphabets is defined on the complete functions:
\begin{equation}
   \sigma_{1}(\mathbb{X} \mathbb{Y}) = \sum_{i} S_{i}(\mathbb{X} \mathbb{Y}) = \prod_{x \in \mathbb{X}} \prod_{y \in \mathbb{Y}} \frac{1}{1-xyt}.
\end{equation}
The function $\sigma_{1}(\mathbb{X} \mathbb{Y})$, denoted by $K(\mathbb{X} , \mathbb{Y})$, is the \emph{Cauchy kernel}. It is the reproducing kernel associated to the usual scalar product on symmetric functions. Therefore, it satisfies the following property: for two basis $A_{\lambda}$ and $B_{\lambda}$ in duality, and two alphabets $\mathbb{X}$ and $\mathbb{Y}$, one has
\begin{equation}
   K(\mathbb{X} , \mathbb{Y}) = \sum_{\lambda} A_{\lambda}(\mathbb{X}) B_{\lambda}(\mathbb{Y}).
\end{equation}
In particular, for Schur functions $S_\lambda=\det\left(S_{\lambda_i-i+j}\right)_{i,j}$, one has
\begin{equation}
\nonumber
K(\X,\Y)=\sum_\lambda S_\lambda(\X)S_\lambda(\Y).   
\end{equation}

\subsection{Jack polynomials}
Jack polynomials $P_{\lambda}^{\left( \frac{1}{\kappa} \right)}$ are defined as the unique family of symmetric functions orthogonal \emph{w.r.t.} a one-parameter deformation $\langle , \rangle_{\frac{1}{\kappa}}$ of the usual scalar product such that, for two partitions $\lambda$ and $\mu$, 
 \begin{equation}
 \label{scalprod}
   \langle p_{\lambda} , p_{\mu} \rangle_{\frac{1}{\kappa}} = z_{\lambda} \left(\frac{1}{\kappa}\right)^{\ell(\lambda)} \delta_{\lambda,\mu},
\end{equation}
(where $\displaystyle{z_{\lambda} = \prod_{i \geq 1} i^{m_{i}(\lambda)} m_{i}(\lambda)!}$ and $\displaystyle{p_{\mu} = \prod_{i=1}^{\ell(\lambda)} p_{\mu_{i}}}$) and such that
\begin{equation}
\label{norma}
 P_{\lambda}^{\left ( \frac{1}{\kappa} \right)}(\mathbb{X}) = m_{\lambda}(\mathbb{X}) + \sum_{\mu \leq \lambda} u_{\lambda,\mu} m_{\mu}(\mathbb{X}).
\end{equation}
where $m_{\lambda}$ denote the monomial symmetric functions.

By definition, it is possible to construct the basis $P_{\lambda}^{\left( \frac{1}{\kappa} \right) }$ by applying Gram-Schmidt algorithm to the Schur basis \emph{w.r.t.} $\langle , \rangle_{\frac{1}{\kappa}}$ and beginning with the partition $1^{n}$. For example, let us construct the $P^{\left(\frac{1}{\kappa}\right)}$ basis of symmetric functions of degree 3. We begin with
\begin{equation}
\nonumber
 P_{111}^{\left(\frac{1}{\kappa}\right)} = S_{111}.
\end{equation}
Then, 
\begin{equation}
\nonumber
 P_{21}^{\left(\frac{1}{\kappa} \right)} = S_{21} - \frac{\langle S_{21} , P_{111}^{\left( \frac{1}{\kappa} \right)} \rangle_{\frac{1}{\kappa}}}{\langle P_{111}^{\left(\frac{1}{\kappa} \right)} , P_{111}^{\left( \frac{1}{\kappa} \right)} \rangle_{\frac{1}{\kappa}}} P_{111}^{\left( \frac{1}{\kappa} \right)}.
\end{equation}
The scalar products $\langle S_{21} , P_{111}^{\left( \frac{1}{\kappa} \right)} \rangle_{ \frac{1}{\kappa} }$ and $\langle P_{111}^{ \left( \frac{1}{\kappa} \right) } , P_{111}^{\left(\frac{1}{\kappa} \right)} \rangle_{\frac{1}{\kappa}}$ are known from the previous steps and the decomposition of Schur functions in the power sums basis. 
The orthogonality property (\ref{scalprod}) of the power sums gives the coefficients.\\
Finally, 
\begin{equation}
\nonumber
 P_{3}^{\left(\frac{1}{\kappa} \right)} = S_{3} - \frac{\langle S_{3} , P_{21}^{\left( \frac{1}{\kappa} \right)} \rangle_{\frac{1}{\kappa}}}{\langle P_{21}^{ \left( \frac{1}{\kappa} \right)} , P_{21}^{ \left( \frac{1}{\kappa} \right)} \rangle_{\frac{1}{\kappa}}} P_{21}^{\left( \frac{1}{\kappa} \right)} - \frac{\langle S_{3} , P_{111}^{\left(\frac{1}{\kappa} \right)} \rangle_{\frac{1}{\kappa}}}{\langle P_{111}^{\left(\frac{1}{\kappa} \right)} , P_{111}^{\left(\frac{1}{\kappa} \right)} \rangle_{\frac{1}{\kappa}}} P_{111}^{\left( \frac{1}{\kappa} \right)}.
\end{equation}
\begin{note}
   Remark that the normalization property (\ref{norma}) is verified since 
\begin{equation}
s_{\lambda} = m_{\lambda} + \sum_{\mu \leq \lambda} w_{\lambda,\mu} m_{\mu}.                                                        
\end{equation}
\end{note}

\subsection{Macdonald polynomials}
\label{Macdo}
The family of symmetric and homogeneous Macdonald $P_{\lambda}(q,t)$ polynomials is uniquely defined by its orthogonality \emph{w.r.t.} the $(q,t)$-deformation $\langle , \rangle_{q,t}$ of the usual scalar pro\-duct on symmetric functions such that
\begin{equation}
   \langle p_{\lambda} , p_{\mu} \rangle_{q,t} = z_{\lambda} \prod_{i=1}^{\ell(\lambda)} \frac{1-q^{\lambda_{i}}}{1-t^{\lambda_{i}}} \delta_{\lambda,\mu},
\end{equation}
and by its dominance property: the expansion of $P_{\lambda}(q,t)$ in terms of $m_{\lambda}$ satisfies
\begin{equation}
 P_{\lambda}(\mathbb{X};q,t) = m_{\lambda}(\mathbb{X}) + \sum_{\mu \leq \lambda} v_{\lambda,\mu} m_{\mu}(\mathbb{X}).
\end{equation}

Denoting by $K_{q,t}(\mathbb{X},\mathbb{Y})$ the reproducing kernel associated to $\langle , \rangle_{q,t}$, one has
\begin{equation}
 K_{q,t}(\mathbb{X},\mathbb{Y}) = \sigma_{1}(\mathbb{X}\mathbb{Y}\frac{1-q}{1-t}),
\end{equation}
and
\begin{equation}
 K_{q,t}(\mathbb{X},\mathbb{Y}) = \sum_{\lambda} P_{\lambda}(\mathbb{X},q,t) Q_{\lambda}(\mathbb{Y},q,t),
\end{equation}
where $Q_{\lambda}(q,t)$ is the dual basis of the $P_{\lambda}(q,t)$ for the scalar product $\langle , \rangle_{q,t}$. Equivalently,
\begin{equation}
   \langle P_{\lambda} (q,t) , Q_{\mu} (q,t) \rangle_{q,t} = \delta_{\lambda,\mu}.
\end{equation}
Note that $Q_\lambda$ and $P_\lambda$ are proportional (by definition).\\
Macdonald polynomials are a generalization of Jack polynomials: substituting $q$ by $t^{\frac{1}{\kappa}}$ and then taking the limit $t \rightarrow 1$, we map $P_{\lambda}(q,t)$ onto $P_{\lambda}^{\left( \frac{1}{\kappa}\right)}$.\\

\section{Exact computation of the Selberg-Jack integral}
\label{exactcomputation}
In this section, we are interested in the computation of 
\begin{equation}
\langle P_{\lambda}^{\left( \frac{1}{\kappa} \right)} \rangle_{a,b,\kappa,N}^{\sharp} : = \frac{\langle P_{\lambda}^{\left( \frac{1}{\kappa} \right)} \rangle_{a,b,\kappa}^{N}}{\langle 1 \rangle_{a,b,\kappa}^{N}}.
\end{equation}
It is clear from (\ref{intMacdo}) that the asymptotics can not be precised directly since the number of factors of the products depends on $N$. Therefore, we have to simplify this expression. 

\begin{lemma}
For any $\kappa \in \mathbb{C}$,
\begin{equation}
\label{basisint}
\begin{aligned}
 \langle P_{\lambda}^{\left( \frac{1}{\kappa} \right)} \rangle_{a,b,\kappa,N}^{\sharp} & = \left[ \prod_{i=1}^{\ell(\lambda)} \prod_{j=i+1}^{\ell(\lambda)} \frac{\Gamma(\lambda_{i}-\lambda_{j}+\kappa(j-i+1)) \Gamma(\kappa(j-i+1))}{\Gamma(\lambda_{i}-\lambda_{j}+\kappa(j-i)) \Gamma(\kappa(j-i))} \right] \\ 
& \left[ \prod_{i=1}^{\ell ( \lambda )} \prod_{j=0}^{\lambda_{i}-1} \frac{\kappa(N+1-i) + j}{\kappa(\ell( \lambda ) + 1 - i) + j} \right]
\left[ \prod_{i=1}^{\ell ( \lambda ) } \prod_{j=0}^{\lambda_{i} - 1} \frac{a + \kappa ( N - i ) + j }{ a + b + \kappa ( 2 N - i - 1 ) + j} \right].
\end{aligned}
\end{equation} 
\end{lemma}
\begin{note} 
It is easy to check that we recover equation (4) of \cite{oai:arXiv.org:1003.5996} by setting $\kappa=1$ in the previous formula.   
\end{note}
\textbf{Proof:} 
We need the following property of the Gamma function:
\begin{equation}
\label{gamma}
 \Gamma(z+n) = \prod_{i=0}^{n-1} (z+i) \Gamma(z).
\end{equation}
First, we split the first product of (\ref{intMacdo}) into several parts: $\displaystyle{\prod_{i=1}^{\ell(\lambda)} \prod_{j=i+1}^{\ell(\lambda)}}$, $\displaystyle{\prod_{i=1}^{\ell(\lambda)} \prod_{j=\ell(\lambda)+1}^{N}}$ and $\displaystyle{\prod_{i=\ell(\lambda)+1}^{N} \prod_{j=i+1}^{N}}$. The first part appears directly in (\ref{basisint}). The last part does not depend on the partition and gives 1 when we divide $\langle P_{\lambda}^{\left( \frac{1}{\kappa} \right)} \rangle_{a,b,\kappa}^{N}$ by $\langle 1 \rangle_{a,b,\kappa}^{N}$. \\
Assuming that $\kappa \in \mathbb{N}$, we can write the second part as
\begin{equation}
\label{inter} 
\prod_{t=0}^{\kappa-1} \prod_{i=1}^{\ell ( \lambda )} \prod_{j=\ell ( \lambda ) +1}^{N} \frac{\lambda_{i}+\kappa(j-i)+t}{\kappa(j-i)+t} .
\end{equation}
Let us define $\displaystyle{W:=\prod_{t=0}^{\kappa-1} \prod_{i=1}^{\ell ( \lambda )} \prod_{j=\ell ( \lambda ) +1}^{N} \frac{\lambda_{i}+\kappa(j-i)+t}{\kappa(j-i)+t}}$. One has
\begin{equation}
\label{trick}
 W = \prod_{i=1}^{\ell ( \lambda )}
\prod_{t=0}^{\lambda_{i}-1} \frac{\kappa(N+1-i) + t}{\kappa(\ell( \lambda ) + 1 - i) + t}.
\end{equation}
To prove the relation, we begin by shifting the variable $t$ in the numerator: $t'=t+\lambda_{i}$:
\begin{equation}
\nonumber
  W = \prod_{i=1}^{\ell(\lambda)} \displaystyle{\frac{\displaystyle{\prod_{t=\lambda_{i}}^{\lambda_{i}+\kappa-1} \prod_{j=\ell(\lambda)+1}^{N} \left( t + \kappa (j-i) \right)}}{\displaystyle{\prod_{t=0}^{\kappa-1} \prod_{j=\ell(\lambda)+1}^{N} \left( t + \kappa (j-i) \right)}}}.
\end{equation}
If $\lambda_{i}>\kappa - 1$, we can multiply numerator and denominator by the factors for $\kappa \leq t \leq \lambda_{i}$: 
\begin{equation}
   W = \prod_{i=1}^{\ell(\lambda)} \prod_{j=\ell(\lambda)+1}^{N} \frac{\displaystyle{\prod_{t=\lambda_{i}}^{\lambda_{i} + \kappa -1} (\kappa (j-i) + t)}}{\displaystyle{\prod_{t=0}^{\kappa -1} (\kappa (j-i) + t)}} \frac{\displaystyle{\prod_{t=\kappa}^{\lambda_{i} -1} (\kappa (j-i) + t)}.}{\displaystyle{\prod_{t=\kappa}^{\lambda_{i} -1} (\kappa (j-i) + t)}}
\end{equation}
If $\lambda_{i}\leq \kappa - 1$, we can also decompose $W$ following the previous equation. In both cases, we have
\begin{equation}
\nonumber
\begin{aligned}
 W = & \prod_{i=1}^{\ell(\lambda)} \displaystyle{\frac{\displaystyle{\prod_{t=\kappa}^{\kappa+\lambda_{i}-1} \prod_{j=\ell(\lambda)+1}^{N} \left( t + \kappa (j-i) \right)}}{\displaystyle{\prod_{t=0}^{\lambda_{i}-1} \prod_{j=\ell(\lambda)+1}^{N} \left( t + \kappa (j-i) \right)}}} \\
= & \prod_{i=1}^{\ell(\lambda)} \displaystyle{\frac{\displaystyle{\prod_{t=0}^{\lambda_{i}-1} \prod_{j=\ell(\lambda)+1}^{N} \left( t + \kappa (j-i+1) \right)}}{\displaystyle{\prod_{t=0}^{\lambda_{i}-1} \prod_{j=\ell(\lambda)+1}^{N} \left( t + \kappa (j-i) \right)}}} .
\end{aligned}
\end{equation}
using another shift of the variable $t$: $t=t'+\kappa.$ A last shift over $j$ gives the following relation:
\begin{equation}
\label{W}
W = \prod_{i=1}^{\ell(\lambda)} \displaystyle{\frac{\displaystyle{\prod_{t=0}^{\lambda_{i}-1} \prod_{j=\ell(\lambda)+2}^{N+1} \left( t + \kappa (j-i) \right)}}{\displaystyle{\prod_{t=0}^{\lambda_{i}-1} \prod_{j=\ell(\lambda)+1}^{N} \left( t + \kappa (j-i) \right)}}}.
\end{equation}
The terms that are different in the numerator and the denominator are obtained for $j=\ell(\lambda)+1$ and $N+1$. After simplification, we obtain (\ref{trick}). \\
To prove (\ref{basisint}), it is enough to use (\ref{gamma}) and the definition of $\langle P_{\lambda}^{\left( \frac{1}{\kappa} \right)} \rangle_{a,b,\kappa,N}^{\sharp}$.
Now, since both left- and right- hand sides of (\ref{W}) are entire functions of exponential type of $N$, 
Carlson's theorem \cite{Carlsson} allows us to extend the result for all $\kappa\in\C$. $\Box$

From (\ref{basisint}) and (\ref{inter}), it is possible to give an approximation of the asymptotic behavior of $\langle P_{\lambda}^{\left( \frac{1}{\kappa} \right)} \rangle_{a,b,\kappa,N}^{\sharp}$ when $N \rightarrow \infty$. Indeed, $\langle P_{\lambda}^{\left( \frac{1}{\kappa} \right)} \rangle_{a,b,\kappa,N}^{\sharp}$ is a product of rational fractions in $N$. Since $N$ appears in the numerator and the denominator of the last product and not at all in the first, the asymptotic behavior depends only on the second product. In this product, $N$ appears to the power $|\lambda|$. Therefore,
\begin{proposition}
\begin{equation}
   \langle P_{\lambda}^{\left( \frac{1}{\kappa} \right)} \rangle_{a,b,\kappa,N}^{\sharp} \underset{N \rightarrow \infty}{\sim} N^{|\lambda|}.
\end{equation}
\end{proposition}

\section{Asymptotics of the Selberg - power sum integral}
\subsection{General method}
In this section, we give the algorithm that allows the computation of $\displaystyle{\langle f \rangle^{\sharp}_{a,b,\kappa,N} = \frac{\langle f \rangle_{a,b,\kappa}^{N}}{\langle 1 \rangle_{a,b,\kappa}^{N}}}$ for any polynomial $f$, not necessarily symmetric. First, we symmetrize $f$:
\begin{equation}
   {\cal S} f = \frac{1}{N!} \sum_{\sigma \in \mathfrak{S}_{n}} \sigma f = \frac{1}{N!} \sum_{\sigma \in \mathfrak{S}_{n}} f(x_{\sigma(1)} , \dots , x_{\sigma(N)}),
\end{equation}
and we remark that $\langle {\cal S} f \rangle^\sharp_{a,b,\kappa,N} =\langle f \rangle^\sharp_{a,b,\kappa,N}$.
The algorithm to compute the integral will then be as follows:
\begin{enumerate}
\item Expand ${\cal S} f$ in terms of Jack $P_{\lambda}^{\left( \frac{1}{\kappa} \right)}$ polynomials;
\item Replace each occurrence of $P_{\lambda}^{\left( \frac{1}{\kappa} \right)}$ by $\langle P_{\lambda}^{\left( \frac{1}{\kappa} \right)} \rangle_{a,b,\kappa,N}^{\sharp}$ which can be computed with (\ref{basisint}).
\end{enumerate}
If $\S f$ admits a decomposition in the Jack basis which does not depend on the number of variables, the asymptotics can be exactly computed since the resulting integral is a rational function in $N$.

We are interested in the value of $\displaystyle{\langle p_{k} \rangle^{\sharp}_{a,b,\kappa,N}}$. In this case, we do not need to symmetrize the integrand and it is sufficient to find the value of the coefficient $\alpha_{\lambda,k}$ of $P_{\lambda}^{\left( \frac{1}{\kappa} \right)}$ in $p_{k}$. This is the aim of the next section and we will find it by manipulating a specialization of Macdonald polynomials. 

\subsection{About a specialization of Macdonald polynomials}
\label{spec}
The aim of this section is to show how to compute the coefficient of the Macdonald polynomial $P_{\lambda}^{q,t}$ in $p_{k}$. 
\begin{proposition}
With the following notations: 
\begin{itemize}
   \item $a_{\lambda}(s)$ denotes the arm-length of $s$ (if $s=(i,j)$, $a_{\lambda}(s) = \lambda_{i}-j$);
   \item $l_{\lambda}(s)$ the leg-length of $s$ (if $s=(i,j)$, $l_{\lambda}(s) = \lambda'_{j}-i$),
\end{itemize}
the coefficient of $P_{\lambda}(q,t)$ in $p_{k}$ is given by
\begin{equation}
   (1-q^{k}) \frac{\displaystyle{\prod_{\genfrac{}{}{0pt}{}{(i,j)\in \lambda}{(i,j)\neq (1,1)}} (t^{i-1}-q^{j-1})}}{ \displaystyle{\prod_{s\in\lambda} (1-q^{a_{\lambda}(s)+1}t^{l_{\lambda}(s)})}}.
\end{equation}
\end{proposition}
\textbf{Proof :} The coefficient of $P_\lambda(\X;q,t)$ in $p_k(\X)$ equals the coefficient of $c_k^{q,t}(\X)$ in $Q_\lambda(\X;q,t)$ where $c_\mu^{q,t}=z_{\mu}(q,t)^{-1}p_\mu(\X)$ denotes the dual basis of $(p_\mu)$ for the scalar product $\langle\,,\,\rangle_{q,t}$ and with
\begin{equation}
   z_{\lambda}(q,t) = z_{\lambda} \prod_{i=1}^{\ell(\lambda)} \frac{1-q^{\lambda_{i}}}{1-t^{\lambda_{i}}}.
\end{equation}
If $\beta_{\lambda,\mu}^{q,t}$ denotes the coefficient of $\displaystyle{p_{\mu}}$ in $Q_\lambda(\X;q,t)$, one has, because of the interpretation of the alphabet $\frac{1-u}{1-t}$ in terms of operations on symmetric functions (see \ref{cauchykern}),
\begin{equation}
Q_\lambda(\frac{1-u}{1-t};q,t)=\sum \beta_{\lambda,\mu}^{q,t} \prod_{i=1}^{\ell(\mu)} \frac{1-u^{\mu_i}}{1-t^{\mu_i}}.
\end{equation}
Dividing left- and right- hand sides of this equation by $\displaystyle{\frac{1}{1-u}}$, we can write that
\begin{equation}
\nonumber
   \frac{1}{1-u}Q_\lambda(\frac{1-u}{1-t};q,t) = \sum_{\mu \leq \lambda} \beta_{\lambda,\mu}^{q,t} \frac{1-u^{\mu_1}}{(1-t^{\mu_1})(1-u)} \prod_{i=2}^{\ell(\mu)} \frac{1-u^{\mu_i}}{1-t^{\mu_i}}.
\end{equation}
The only term that does not vanish for $u=1$ is $\displaystyle{\frac{1-u^{\mu_1}}{(1-t^{\mu_1})(1-u)}}$ because of the pole at this point. Therefore, if we take the limit $u \rightarrow 1$, the partitions that give a non zero contribution have one part: $\mu=|\lambda|=k$. Since 
\begin{equation}
\nonumber
   \lim_{u \rightarrow 1} \frac{1-u^{k}}{1-u}=k,
\end{equation}
the following relation holds:
\begin{equation}
\nonumber
\frac{k}{1-t^k} \beta_{\lambda,k}^{q,t}=\lim_{u\rightarrow1}\frac{1}{1-u}Q_\lambda(\frac{1-u}{1-t};q,t)   .
\end{equation}
But $Q_\lambda(\frac{1-u}{1-t};q,t)$ is known and can be expressed with help of formulae (8.3) p352 and (8.8) p 354 of \cite{Macdonald95}:
\begin{equation}
\nonumber
Q_\lambda(\frac{1-u}{1-t};q,t)=\frac{\displaystyle{\prod_{(i,j)\in\lambda}(t^{i-1}-q^{j-1}u)}}{\displaystyle{\prod_{s\in\lambda}
(1-q^{a_{\lambda}(s)+1}t^{l_{\lambda}(s)})}}.
\end{equation}
Therefore, 
\begin{equation}
   \beta_{\lambda,k}^{q,t} = \frac{1-t^k}{k} \frac{\displaystyle{\prod_{\genfrac{}{}{0pt}{}{(i,j)\in \lambda}{(i,j)\neq (1,1)}} (t^{i-1}-q^{j-1})}}{ \displaystyle{\prod_{s\in\lambda} (1-q^{a_{\lambda}(s)+1}t^{l_{\lambda}(s)})}},
\end{equation}
and the coefficient of $P_\lambda(\X;q,t)$ in $p_k(\X)$ is equal to 
\begin{equation}
\begin{aligned}
z_k(q,t) \beta_{\lambda,k}^{q,t} & = k \frac{1-q^{k}}{1-t^{k}} \frac{1-t^k}{k} \frac{\displaystyle{\prod_{\genfrac{}{}{0pt}{}{(i,j)\in \lambda}{(i,j)\neq (1,1)}} (t^{i-1}-q^{j-1})}}{ \displaystyle{\prod_{s\in\lambda} (1-q^{a_{\lambda}(s)+1}t^{l_{\lambda}(s)})}} \\
 & = (1-q^{k}) \frac{\displaystyle{\prod_{\genfrac{}{}{0pt}{}{(i,j)\in \lambda}{(i,j)\neq (1,1)}} (t^{i-1}-q^{j-1})}}{ \displaystyle{\prod_{s\in\lambda} (1-q^{a_{\lambda}(s)+1}t^{l_{\lambda}(s)})}}\Box.
\end{aligned}
\end{equation} 
\subsection{Exact computation}
The result of the previous section involves Macdonald polynomials. We have to come back to Jack polynomials. Applying the method described in section \ref{Macdo} (set $\displaystyle{q=t^{\frac{1}{\kappa}}}$ and take the limit $t \rightarrow 1$), one recovers a result of \cite{StanleyStembridge}:
\begin{corollary}
The coefficient of $P_{\lambda}^{\left( \frac{1}{\kappa} \right)}$ in $p_{k}$ is equal to 
\begin{equation}
\alpha_{\lambda,k} = k \frac{\displaystyle{\prod_{\genfrac{}{}{0pt}{}{(i,j)\in \lambda}{(i,j)\neq (1,1)}} ((j-1)-\kappa(i-1))}}{ \displaystyle{\prod_{s\in\lambda} \left( a_{\lambda}(s)+1 + l_{\lambda}(s) \kappa \right)}}.
\end{equation}
\end{corollary}


Applying the algorithm, we have the following proposition:
\begin{proposition}
\begin{equation}
   \langle p_{k} \rangle_{a,b,\kappa,N}^{\sharp} = k \sum_{\lambda \vdash k} \frac{\displaystyle{\prod_{\genfrac{}{}{0pt}{}{(i,j)\in \lambda}{(i,j)\neq (1,1)}} ((j-1)-\kappa(i-1))}}{ \displaystyle{\prod_{s\in\lambda} \left( a_{\lambda}(s)+1 + l_{\lambda}(s) \kappa \right)}} \langle P_{\lambda}^{\left( \frac{1}{\kappa} \right)} \rangle_{a,b,\kappa,N}^{\sharp}.
\end{equation}
\end{proposition}

We have written the integral $\langle p_k \rangle^\sharp_{a,b,\kappa,N}$ as a rational fraction in $N$. Now, the study of the asymptotics is easy whence simplifying this fraction. Numerical evidences suggest that, if $a$, $b$ and $\kappa$ are constant, $\displaystyle{\{p_k\}_{a,b,\kappa}:=\lim_{N \rightarrow \infty} \frac{1}{N} \langle p_k \rangle^\sharp_{a,b,\kappa,N}}$ does not depend on $a$, $b$ and $\kappa$ and its value involves central binomial coefficients:
\begin{equation}
\{p_k\}_{a,b,\kappa}=\frac{1}{2^{2k}}\genfrac{(}{)}{0pt}{}{2k}{k}.   
\end{equation}

Even when $a$ and $b$ are linear in $N$, the limit of $\frac{1}{N} \langle p_k \rangle_{a,b,\kappa,N}^{\sharp}$ seems to converge. In particular, when $a=\kappa(\ell-1)N$ and $b=0$, we find experimentally that its value is related to the number of symmetric Dyck paths counted by number of peaks:
\begin{equation}
\lim_{N \rightarrow \infty} \frac{1}{N} \langle p_k \rangle^\sharp_{\kappa(\ell-1)N,0,\kappa,N} = \frac{\ell}{(1+\ell)^{2k-1}} \sum_{i=0}^{2(k-1)} \genfrac{(}{)}{0pt}{}{k-1}{\left\lceil \frac{i}{2} \right\rceil}
\genfrac{(}{)}{0pt}{}{k-1}{\left\lfloor \frac{i}{2} \right\rfloor}
\ell^i.
\end{equation}

In \cite{oai:arXiv.org:1003.5996}, we have completely described, in a combinatorial way, what happens when $\kappa=1$ and, in a forthcoming paper, we will prove the result for formal parameter $\kappa$ using analytic tools.

\bibliographystyle{amsplain}
\bibliography{ASelberg.bib}

\end{document}